\documentclass[11pt]{article}
\usepackage{amsfonts,amssymb,amsmath}
\newtheorem{definition}{Definition}
\newtheorem{proposition}[definition]{Proposition}

\DeclareMathOperator{\rank}{rank}
\DeclareMathOperator{\divz}{div}
\DeclareMathOperator{\curl}{curl}

\begin{document}

\title {Ruled surfaces asymptotically normalized}
\author {Stylianos Stamatakis, Ioannis Kaffas\\
            Department of Mathematics, Aristotle University of Thessaloniki\\
                GR-54124 Thessaloniki, Greece\\
                e-mail: stamata@math.auth.gr}
\date{}
\maketitle

\begin{abstract}
                We consider a skew ruled surface $\Phi$ in the Euclidean space $E^{3}$
                and relative normalizations of it, so that the relative normals at each
                point lie in the corresponding asymptotic plane of $\Phi$. We call such
                relative normalizations and the resulting relative images of $\Phi$ \emph{asymptotic}.
                We determine all ruled surfaces and the asymptotic normalizations of them,
                for which $\Phi$ is a relative sphere (proper or inproper) or the asymptotic
                image degenerates into a curve. Moreover we study the sequence of the ruled
                surfaces $\left \{\Psi_{i}\right\}_{i\in \mathbb{N}}$, where $\Psi_{1}$ is an
                asymptotic image of $\Phi$ and $\Psi_{i}$, for $i\geq2$, is an asymptotic image
                of $\Psi_{i-1}$. We conclude the paper by the study of various properties
                concerning some vector fields, which are related with $\Phi$.
                \\[1mm] {\em MSC 2010:} 53A25,\ 53A05, 53A15, 53A40
                \\[1mm] {\em Keywords:} Ruled surfaces, relative normalizations
\end{abstract}

\section{Preliminaries}

Here we sum up briefly some elementary facts concerning the relative Differential Geometry
of surfaces and the Differential Geometry of ruled surfaces in the Euclidean space $E^{3}$;
for notations and definitions the reader is referred to \cite{Pottmann} and \cite{Schirokow}.


In the Euclidean space $E^{3}$ let $\Phi:\bar{x}=\bar{x}(u,v)$ be an injective $C^{r} $-immersion defined on a
region $U $ of $\mathbb{R}^{2}$, with non-vanishing Gaussian curvature. A $C^{s}$-mapping $\bar{y}:U\longrightarrow E^{3},\, r>
s\geq1$, is called a $C^{s}$-relative normalization of $\Phi$ if
\begin{equation}
\rank \left(\{\bar{x}_{/1},\bar{x}_{/2},\bar{y}\}\right)=3,\,\,
\rank \left(\{\bar{x}_{/1},\bar{x}_{/2},\bar{y}_{/i}\}\right)=2,\,\,
i=1,2,\,\,\forall \left(u,v\right) \in U,  \label{1}
\end{equation}%
where
$$
f_{/i}:=\frac{\partial f}{\partial u^{i}},~f_{/ij}:=\frac{\partial ^{2}f}{\partial u^{i}\partial u^{j}} \quad \text{etc.}
$$
denote partial derivatives of a function (or a vector-valued function) $f$ in the coordinates
$u^{1}:=u,\,u^{2}:=v$. The covector $\bar{X} $ of the tangent plane is defined by
\begin{equation}
\langle \bar{X},\bar{x}_{/i}\rangle =0 \quad (i=1,2)\quad \text{and}\quad
\langle \bar{X},\bar{y}\rangle =1,  \label{2}
\end{equation}%
where $\langle \,,\,\rangle$ denotes the standard scalar product in $E^{3}$. The relative metric $G$ is introduced by%
\begin{equation}
G_{ij}=\langle \bar{X},\bar{x}_{/ij}\rangle.  \label{3}
\end{equation}%
The \emph{support function of the relative normalization} $\bar{y} $ is defined by $q:=\langle \bar{\xi},\bar{y}\rangle$ (see \cite{Manhart3}), where $\bar{\xi} $ is the Euclidean normalization of $\Phi $. By virtue of (\ref{1}) $q $ never vanishes on $U $ and, because of (\ref{2}),
$\bar{X}=q^{-1}\bar{\xi} $. Then by (\ref{3}), we also obtain%
\begin{equation}
G_{ij}=q^{-1}\,h_{ij},  \label{12}
\end{equation}%
where $h_{ij} $ are the coefficients of the second fundamental form of $\Phi$. Conversely, when a
support function $q $ is given, then the relative normalization $\bar{y} $ is uniquely determined by (see
\cite[p. 197]{Manhart3})
\begin{equation}
\bar{y}=-h^{\left( ij\right) }\,q_{/i}\,\,\bar{x}_{/j}+q\,\bar{\xi},  \label{6}
\end{equation}%
where $h^{(ij)} $ are the coefficients of the inverse tensor of $h_{ij} $. For a function (or a vector-
valued function) $f $ we denote by $\nabla ^{G}f$ the first Beltrami differential operator and by $\nabla _{i}^{G}f $ the
covariant derivative, both with respect to the relative metric. We consider the coefficients
$$
A_{ijk}:=\langle \bar{X},\nabla _{k}^{G}\nabla _{j}^{G}\bar{x}_{/i}\rangle
$$%
of the Darboux tensor. Then, by using the relative metric tensor $G_{ij} $ for ``raising and lowering
the indices'', the Tchebychev vector $\bar{T}$ of the relative normalization $\bar{y}$ is defined by%
\begin{equation}
\bar{T}:=T^{m}\, \bar{x}_{/m}\quad \text{where\quad }T^{m}:=\frac{1}{2}%
A_{i}^{im} \label{10}
\end{equation}%
and the Pick invariant by%
\begin{equation}
J:=\frac{1}{2}A_{ijk}A^{ijk}.  \label{26}
\end{equation}%
The relative shape operator has the coefficients $B_{i}^{j}$ defined by
\begin{equation}
\bar{y}_{/i}=:-B_{i}^{j}\, \bar{x}_{/j}.  \label{27}
\end{equation}%
Then, the relative curvature and the relative mean curvature are defined by
\begin{equation}
K:=\det \left(B_{i}^{j}\right),\quad H:=\frac{B_{1}^{1}+B_{2}^{2}}{2}.  \label{28}
\end{equation}%
When we attach the vectors $\bar{y} $ of the relative normalization to the origin, the endpoints of
them describe the relative image of $\Phi$.

\vspace{2mm}
Let now $\Phi $ be a skew (non-developable) ruled $C^{2} $-surface, which is defined by its striction
curve $\Gamma: \bar{s}=\bar{s}(u),\,u\in I \,\,(I\subset \mathbb{R} $ open interval) and the unit vector $\bar{e}$ pointing along the generators. We choose the parameter $u $ to be the arc length along the spherical curve $\bar{e}=\bar{e}(u)$
and we denote the differentiation with respect to $u $ by a prime. Then a parametrization of
the ruled surface $\Phi $ over the region $U:=I\times \mathbb{R} $ is
\begin{equation}
\bar{x}(u,v)=\bar{s}(u)+v\,\bar{e}(u),  \label{4}
\end{equation}%
with%
\begin{equation}
\left \vert \bar{e}\right \vert =|\bar{e}^{\prime }|=1,\quad \langle \bar{s}%
^{\prime }(u),\bar{e}^{\prime }(u)\rangle =0\quad \forall ~u\in I.
\label{23}
\end{equation}%
The distribution parameter $\delta(u) :=(\bar{s}^{\prime },\bar{e},\bar{e}^{\prime })$, the conical curvature $\kappa(u):=(\bar{e},\bar{e}^{\prime},
\bar{e}^{\prime \prime })$ and the
function $\lambda :=\cot \sigma$, where $\sigma :=\sphericalangle (\bar{e},\bar{s}^{\prime })$ is the striction of $\Phi\,\,(-\frac{\pi }{2}<\sigma \leq \frac{\pi }{2}, \,\,{\textrm{sign}}\sigma =\textrm{sign}\delta)$,
are the fundamental invariants of $\Phi$ and determine uniquely, up to Euclidean rigid motions,
the ruled surface $\Phi $. The moving frame of $\Phi $ is the orthonormal frame which is attached to
the striction point $\bar{s}(u)$, and consists of the vector $\bar{e}(u) $, the central normal vector $\bar{n}(u):=\bar{e}^{\prime }(u) $
and the central tangent vector $\bar{z}(u):=\bar{e}(u)\times \bar{n}(u) $. It fulfils the equations \cite[p. 280]{Pottmann}
\begin{equation}
\bar{e}^{\prime }=\bar{n},\quad \bar{n}^{\prime }=-\bar{e}+\kappa \,\bar{z},\quad \bar{z}^{\prime }=-\kappa \,\bar{n}.  \label{21}
\end{equation}
Then, we have
\begin{equation}
\bar{s}^{\prime }=\delta \lambda \bar{e}+\delta \bar{z}.  \label{5}
\end{equation}%
By (\ref{4}) and (\ref{5})\ we also obtain%
\begin{equation}
\bar{x}_{/1}=\delta \lambda \bar{e}+v\bar{n}+\delta \bar{z},\quad \bar{x}%
_{/2}=\bar{e},  \label{9}
\end{equation}%
and thus%
\begin{equation}
\bar{\xi}=\frac{\delta \bar{n}-v\bar{z}}{w},\quad \text{where}\quad w:=\sqrt{%
v^{2}+\delta ^{2}}.  \label{18}
\end{equation}%
The coefficients $g_{ij} $ and $h_{ij} $ of the first and the second fundamental form of $\Phi $ take the form%
\begin{equation}
g_{11}=w^{2}+\delta ^{2}\lambda ^{2},\quad
g_{12}=\delta \lambda,\quad
g_{22}=1,  \label{47}
\end{equation}%
\begin{equation}
h_{11}=-\frac{\kappa w^{2}+\delta ' v-\delta ^{2}\lambda }{w},\quad
h_{12}=\frac{\delta }{w}, \quad
h_{22}=0.  \label{7}
\end{equation}%
The Gaussian curvature $\widetilde{K} $ of $\Phi $ is given by (E. Larmarle's formula \cite{Pottmann})%
\begin{equation}
\widetilde{K}=-\frac{\delta ^{2}}{w^{4}}.  \label{17}
\end{equation}%
In this paper only skew ruled surfaces of the space $E^{3} $ are considered with parametrization
like in (\ref{4}) and (\ref{23}).

\section{Ruled surfaces relatively normalized}

Let $\bar{y}$ be a relative normalization of a given ruled $C^{2}$-surface $\Phi\,\, (\delta \neq 0$) and let $q$ be the
corresponding support function. Then, on account of (\ref{12}) and (\ref{7}) the coefficients of the
inverse relative metric tensor are computed by
\begin{equation}
G^{(11)}=0,\quad G^{(12)}=\frac{wq}{\delta },\quad G^{(22)}=\frac{wq\left( \kappa
w^{2}+\delta ^{\prime }v-\delta ^{2}\lambda \right) }{\delta ^{2}}.
\label{16}
\end{equation}%

The relative normalization $\bar{y}$ of $\Phi $ can be expressed with respect to the moving frame $\left \{ \bar{e},\bar{n},\bar{z}\right \}$,
by using (\ref{6}), (\ref{9}), (\ref{18}) and (\ref{7}), as follows:
\begin{equation}
\bar{y}=-w\frac{\delta q_{/1}+q_{/2}(\kappa w^{2}+\delta ^{\prime }v)}{%
\delta ^{2}}\bar{e}+\frac{\delta ^{2}q-w^{2}vq_{/2}}{\delta w}\bar{n}-\frac{%
vq+w^{2}q_{/2}}{w}\bar{z}.  \label{8}
\end{equation}%

It is well known \cite[p. 199]{Manhart3}, that the components of the Tchebychev vector $\bar{T}$ of $\bar{y}$ are given by%
\begin{equation}\label{61}
T^{i}=\left[ \ln \left( \frac{|q|}{q_{AFF}}\right) \right] _{/j}G^{(ij)},
\end{equation}%
where, by virtue of (\ref{17}),
\begin{equation}\label{25}
q_{AFF}=|\widetilde{K}|^{1/4}=\frac{|\delta|^{1/2}}{w}
\end{equation}
denotes the support function of the equiaffine normalization $\bar{y}_{AFF}$.
From the relations (\ref{17}) and (\ref{16}) we have
\begin{equation}
T^{1}=\frac{w^{2}q_{/2}+vq}{\delta w}, \quad
T^{2}=\frac{2\delta w^{2}q_{/1}+\delta
^{\prime }q\left( \delta ^{2}-v^{2}\right) }{2\delta ^{2}w}+\frac{%
T^{1}\left( \kappa w^{2}+\delta ^{\prime }v-\delta ^{2}\lambda \right) }{%
\delta }.  \label{43}
\end{equation}%
Thus, by using (\ref{10}) and (\ref{9}), we obtain%
\begin{equation}
\bar{T}=w\frac{q\left( 2\kappa v+\delta ^{\prime }\right) +2\delta
q_{/1}+2q_{/2}(\kappa w^{2}+\delta ^{\prime }v)}{2\delta ^{2}}\bar{e}+\frac{%
vq+w^{2}q_{/2}}{\delta w}\left( v\bar{n}+\delta \bar{z}\right) .  \label{13}
\end{equation}%
Especially, the Tchebychev vector $\bar{T}_{EUK}$ of the Euclidean normalization ($q=1$) reads%
\begin{equation}
\bar{T}_{EUK}=w\frac{2\kappa v+\delta ^{\prime }}{2\delta ^{2}}\bar{e}+\frac{%
v}{\delta w}\left( v\bar{n}+\delta \bar{z}\right) .  \label{14}
\end{equation}

We introduce now the tangential vector
\begin{equation}
\bar{Q}:=\frac{1}{4}\nabla ^{G}\left( \frac{1}{q},\bar{x}\right)  \label{11}
\end{equation}%
of $\Phi$. On account of (\ref{6}) and (\ref{16}) we have
$$
\bar{y}-q\,\bar{\xi}=4\,q\,\bar{Q}.
$$
Thus, by (\ref{11}), the vector $\bar{Q}$ \emph{is in the direction of the tangential component of $\bar{y}$}.

\begin{definition}
We call $\bar{Q}$ the support vector of $\bar{y}$.
\end{definition}

Its components with respect to the local basis \{$\bar{x}_{/1},\bar{x}_{/2}$\}, because of (\ref{16}) and (\ref{11}), are
\begin{equation}
Q^{1}=\frac{-w\,q_{/2}}{4\delta q},\quad Q^{2}=-w\frac{\left( \kappa
w^{2}+\delta ^{\prime }v-\delta ^{2}\lambda \right) q_{/2}+\delta q_{/1}}{%
4\delta ^{2}q}.  \label{45}
\end{equation}%
By using (\ref{9}) we find%
\begin{equation}
\bar{Q}=-w\frac{\delta \, q_{/1}+q_{/2}(\kappa w^{2}+\delta ^{\prime }v)}{%
4\delta ^{2}q}\bar{e}-\frac{wq_{/2}}{4\delta q}\left( v\bar{n}+\delta \bar{z}%
\right).  \label{15}
\end{equation}%
Denoting by $\bar{Q}_{AFF}$ the support vector of the equiaffine normalization $\bar{y}_{AFF}$ and using (\ref{25}), (\ref{13}), (\ref{14}) and (\ref{15}), we get the relations%
$$
\bar{T}_{EUK}=4\bar{Q}_{AFF},\quad \bar{T}=q\bar{T}_{EUK}-4q\bar{Q}.
$$

\section{Asymptotic normalizations of ruled surfaces}

First in this section we find all relative normalizations $\bar{y} $, so that the relative normals at
each point $P $ of $\Phi $ lie in the corresponding asymptotic plane, i.e. in the plane \{$ P;\bar{e},\bar{n}$\}. On
account of (\ref{8}), this is valid iff $vq + w^{2}q_{/2}=0 $, or, equivalently, iff the support function $q $ of
$\bar{y} $ is of the form $ q=f\,w^{-1}$, where $f=f(u) $ is an arbitrary non-vanishing $C^{1} $-function. By
virtue of (\ref{13}) we have

\begin{proposition}
The following statements are equivalent: (a) The relative normals at each
point $P$ of $\Phi$ lie on the corresponding asymptotic plane. (b) The Tchebychev vector $\bar{T}$ of $\bar{y}$
at each point $P$ of $\Phi$ is parallel to the corresponding generator. (c) The support function is
of the form
\begin{equation}
q=\frac{f\left( u\right) }{w},\quad f\left( u\right) \in C^{1}\left(
I\right) ,\quad f\left( u\right) \neq 0.  \label{19}
\end{equation}
\end{proposition}

\begin{definition}
We call a support function of the form (\ref{19}), as well as the corresponding relative normalization
\begin{equation}
\bar{y}=\left[ -\left( \frac{f}{\delta }\right) ^{\prime }+\frac{\kappa f}{%
\delta ^{2}}v\right] \bar{e}+\frac{f}{\delta }\bar{n},  \label{20}
\end{equation}%
and the resulting relative image of $\Phi$ asymptotic.
\end{definition}

It is apparent from (\ref{25}) and (\ref{19}), that the \emph{equiaffine normalization} $\bar{y}_{AFF}$ \emph{is contained in the set of the asymptotic ones}. Support functions of ruled surfaces of the form (\ref{19}) were introduced by the first author in \cite{Stamatakis}.

\vspace{2mm}
We consider an asymptotically normalized by (\ref{20}) ruled surface $\Phi $. The Pick invariant
of $\Phi $ is computed from (\ref{26}), by using the well known equation \cite[p. 196]{Manhart3}%
\begin{equation}
A_{ijk}=\frac{1}{q}\langle \bar{\xi},\bar{x}_{/ijk}\rangle -\frac{1}{2}%
\left( G_{ij/k}+G_{jk/i}+G_{ki/j}\right)  \label{31}
\end{equation}%
and the relations (\ref{12}), (\ref{9}), (\ref{18}) and (\ref{7}). We easily find $A_{222}=0 $. Then, since the Darboux
tensor is fully symmetric, we have%
\begin{equation}
J=\frac{3}{2}\left( A_{112}A^{112}+A_{122}A^{122}\right) .  \label{46}
\end{equation}%
On account of (\ref{31}), by straightforward calculations, we get
$$
A_{112}=\frac{2\delta f^{\prime }-\delta ^{\prime }f}{2f^{2}},\quad
A^{112}=A_{122}=0,\quad A^{122}=f\frac{2\delta f^{\prime }-\delta ^{\prime }f%
}{2\delta ^{3}}.
$$%
Substitution in (\ref{46}) gives $ J=0 $. This generalizes a result
on equiaffinelly normalized ruled surfaces (see \cite[p. 217]{Blaschke}).

\vspace{2mm}
The relative curvature and the relative mean curvature of $\Phi $ are computed on account of
(\ref{28}). By using (\ref{27}), (\ref{9}) and (\ref{20}), we find the coefficients of the relative shape operator
\begin{equation}
B_{1}^{1}=\frac{-\kappa f}{\delta ^{2}},\quad B_{2}^{1}=0,\quad B_{2}^{2}=\frac{-\kappa f}{\delta ^{2}}, \label{70}
\end{equation}
\begin{equation}
B_{1}^{2}=\frac{2\delta ^{\prime }f\left( \kappa v +\delta ^{\prime}\right) -\delta \left[ \kappa f^{\prime }v+
2\delta^{\prime }f^{\prime}+f\left( \kappa^{\prime }v+\delta ^{\prime \prime }\right) \right] +\delta
^{2} \left[ f\left( 1+\kappa \lambda \right) +f^{\prime \prime }\right]}{\delta ^{3}}, \label{71}
\end{equation}%
so that%
\begin{equation}
K=\frac{\kappa ^{2}f^{2}}{\delta ^{4}},\quad H=\frac{-\kappa f}{\delta ^{2}}. \label{59}
\end{equation}%
It is obvious that:
\vspace{2mm}
\begin{itemize}
   \item
            \emph{The relative curvature and the relative mean curvature are constant along each generator of $\Phi$.
            Moreover they are both constant iff the function $f$ is of the form $f{\small{=}}~c\, \delta^{2}\,\kappa^{-1},
            \,c\in ~\mathbb{R}^{\ast}$}.
            \item
            \emph{The only asymptotically normalized ruled surfaces, which are relative minimal surfaces} (\emph{or of vanishing relative curvature}) \emph{are the conoidal ones}.
\end{itemize}

\vspace{2mm}
The scalar curvature $S$ of the relative metric $G$, which is defined formally and is the
curvature of the pseudo-Riemannian manifold ($\Phi,G$), is obtained by direct computation to
be $S=H$. Substituting $J,H $ and $S$ in the \emph{Theorema Egregium of the relative
Differential Geometry} (see \cite[p. 197]{Manhart3}), which states that
$$
H-S+J=2\,T_{i}\,T^{i},
$$%
it turns out that \emph{the norm $\left \Vert T\right \Vert _{G} $ with respect to the relative metric of the Tchebychev vector
$\bar{T}$ of any asymptotic normalization $\bar{y}$ of $\Phi$ vanishes identically}.

\vspace{2mm}
Let the ruled surface $\Phi$ be non-conoidal. We consider the covariant coefficients $B_{ij}{\small{=}}~
B_{i}^{k}\,G_{kj}$ of the relative shape operator and we denote by $\widetilde{B}$ the scalar curvature of the metric
$B{ij}\,du^{i}\,du^{j}$, which is defined formally just as the curvature $S$. Then, on account of (\ref{12}), (\ref{7}),
(\ref{19}), (\ref{70}) and (\ref{71}), it turns out that $\widetilde{B}$ equals 1.

\vspace{2mm}
From (\ref{20}) it is obvious, that the asymptotic image of $\Phi$ degenerates into a point or into
a curve iff $\Phi$ is conoidal. In this case we have
$$
\bar{y}=-\left( \frac{f}{\delta }\right) ^{\prime }\bar{e}+\frac{f}{\delta }%
\bar{n}.
$$%
Furthermore, computing the derivative of $\bar{y} $ and using (\ref{21}), it follows immediately that the
asymptotic image of $\Phi $ degenerates a) into a curve $\Gamma_{1}$, iff $f\neq \delta(c_{1}\cos u + c_{2}\sin u)$, $c_{1},c_{2}\in \mathbb{R}$,
\,$c_{1}^2+c_{2}^2\neq 0$, or b) into a point, iff $f=\delta(c_{1}\cos u+c_{2}\sin u)$, $c_{1},c_{2}\in \mathbb{R}$, $c_{1}^2+c_{2}^2\neq 0$. In case (a) one readily verifies, that $\Gamma _{1} $ is a planar curve, whose radius of curvature equals $r=|(\frac{f}{\delta})^{\prime \prime }+\frac{f}{\delta}|$.
In case (b) the asymptotic normalization of $\Phi$ is constant. Consequently the ruled surface $\Phi $ is an
improper relative sphere \cite{Manhart1}. Hence we have

\begin{proposition}
Let $\Phi$ be an asymptotically normalized ruled surface. The asymptotic image
of $\Phi$ degenerates (a) into a curve, which is planar, iff $\Phi$ is conoidal and $f\neq\delta(c_{1}\cos u+
c_{2}\sin u),\,c_{1},c_{2}\in \mathbb{R},\,\, c_{1}^2+c_{2}^2\neq 0$, (b) into a point, whereupon $\Phi$ is an improper relative sphere,
iff $\Phi $ is conoidal and $f=\delta (c_{1}\cos u+c_{2}\sin u),\,c_{1},c_{2}\in \mathbb{R},\,c_{1}^2+c_{2}^2\neq 0$.
\end{proposition}

Let now $\Phi $ be a proper relative sphere, i.e. its relative normals pass through a fixed point \cite{Manhart2}.
It is well known, that this is valid iff there exists a constant $c\in \mathbb{R}^{\ast}$ and a constant vector
$\bar{a} $, such that $\bar{x}=c\,\bar{y}+\bar{a} $. Taking partial derivatives of this last equation on account of (\ref{4}),
(\ref{21}), (\ref{5}), (\ref{20}) and (\ref{59}), we obtain
\begin{equation}
f=\frac{\delta^2}{c\kappa},\quad (\kappa \neq 0),  \label{22}
\end{equation}%
\begin{equation}
\left( \frac{\delta }{\kappa }\right)^{\prime \prime }+\frac{\delta }{\kappa
}\left( 1+\kappa \lambda \right) =0  \label{32}
\end{equation}%
and%
\begin{equation}
c\bar{y}=\left[ -\left( \frac{\delta }{\kappa }\right)^{\prime }+v\right]
\bar{e}+\frac{\delta }{\kappa }\bar{n}.  \label{24}
\end{equation}%
We notice, that the relative curvature and the relative mean curvature of a proper relative sphere are constant.\\
Conversely, let us suppose, that the equations (\ref{22}) and (\ref{32}) are
valid, where $c\in \mathbb{R}^{\ast}$. Then, because of (\ref{20}), the equation (\ref{24}) is valid as
well. Moreover, from (\ref{5}) and (\ref{32}) we obtain%
$$
\left[ -\left( \frac{\delta }{\kappa }\right)^{\prime }\bar{e}+\frac{\delta
}{\kappa }\bar{n}\right]^{\prime }=\bar{s}^{\prime }.
$$%
Therefore the striction curve $\Gamma $ of $\Phi $ is parametrized by
\begin{equation}
\bar{s}=-\left( \frac{\delta }{\kappa }\right)^{\prime }\bar{e}+\frac{\delta
}{\kappa }\bar{n}+\bar{a},\quad \bar{a}=const.  \label{34}
\end{equation}%
By combining this last relation with (\ref{4}) and (\ref{24}) we get $\bar{x}%
=c\,\bar{y}+\bar{a} $, which means that $\Phi $ is a proper relative sphere.
Thus, we arrive at
\begin{proposition}{\label{Prop}}
An asymptotically normalized ruled surface $\Phi $ is a proper relative sphere iff the function
$f$ is given by (\ref{22}) and its fundamental invariants are related as in the equation (\ref{32}).
\end{proposition}

We now assume, that the relative normals of $\Phi $ are parallel to a fixed
plane $E $. Let $\bar{c} $ be a constant normal unit vector on $E $. Then $%
\langle \bar{y},\bar{c}\rangle =0, $ whence, on account of (\ref{20}), we find%
\begin{equation}
\frac{\kappa f}{\delta ^{2}}\langle \bar{e},\bar{c}\rangle v+\left[ -\left(
\frac{f}{\delta }\right)^{\prime }\langle \bar{e},\bar{c}\rangle +\frac{f}{%
\delta }\langle \bar{n},\bar{c}\rangle \right] =0.  \label{51}
\end{equation}%
Differentiation of (\ref{51}) relative to $v $ yields $\kappa \langle \bar{e}%
,\bar{c}\rangle =0 $. Then, again from (\ref{51}), we derive the system%
$$
\kappa \langle \bar{e},\bar{c}\rangle =0,\quad \left( \frac{f}{\delta }%
\right)^{\prime}\langle \bar{e},\bar{c}\rangle -\frac{f}{\delta }%
\langle \bar{n},\bar{c}\rangle =0.
$$
In case of $\langle \bar{e},\bar{c}\rangle \neq 0$, we obtain%
$$
\kappa =0,\quad \left( \frac{f}{\delta }\right)^{\prime \prime }+\frac{f}{\delta }=0.
$$%
In this case $\bar{y}$ is constant, i.e. $\Phi$ is an improper relative sphere. In case of $\langle \bar{e},\bar{c}\rangle =0 $, we have $\kappa=0$
and (\ref{51}) is identically fulfilled. So we have proved
\begin{proposition}
If the relative normals of an asymptotically normalized ruled surface $\Phi $ are parallel to a fixed plane $E$, then $\Phi$ is conoidal. Furthermore $\Phi$ is either an improper relative sphere or its generators are parallel to $E$.
\end{proposition}

We consider now a non-conoidal ruled surface which is asymptotically normalized by (\ref{20}).
In view of (\ref{59}) we observe that \emph{all points of $\Phi $ are relative umbilics} ($H^{2}-K=0 $), result which generalizes a result on equiaffinelly normalized ruled surfaces (see \cite[p. 218]{Blaschke})
Thus, the relative principal curvatures $k_{1}$ and $k_{2}$ equal $H$. The parametrization of the unique relative focal
surface of $\Phi$, which initially reads%
$$
\bar{x}^{\ast }=\bar{s}+v\bar{e}+\frac{1}{H}\bar{y},
$$%
becomes%
$$
\bar{x}^{\ast }=\bar{s}-\frac{\delta }{\kappa }\bar{n}+\frac{\delta
f^{\prime }-\delta ^{\prime }f}{\kappa f}\bar{e},
$$%
i.e. \emph{the focal surface degenerates into a curve $\Gamma ^{\ast } $ and all relative normals
along each generator form a pencil of straight lines}. This generalizes a result
on equiaffinelly normalized ruled surfaces (see \cite[p. 204]{Schirokow}).\\
Let $P(u_{0})$ be a point of the striction line $\Gamma$ of $\Phi$ and $R(u_{0})$ the corresponding
point on the focal curve $\Gamma^{\ast}$. If we consider all asymptotic normalizations of $\Phi$,
then the locus of the points $R( u_{0})$ is a straight line parallel to the vector $\bar{e}\left( u_{0}\right) $.
In this way we obtain a ruled surface $\Phi ^{\ast }$, whose generators are parallel to the
vectors $\bar{e}(u)$, a parametrization of which reads%
$$
\Phi ^{\ast }:x^{\ast }=\bar{s}-\frac{\delta }{\kappa }\bar{n}+v^{\ast }\,\bar{e},
$$
which is \emph{the asymptotic developable of} $\Phi$ (see \cite[p. 51]{Hoschek}). One easily verifies, that
$$
\bar{s}^{\ast }=\bar{s}-\frac{\delta }{\kappa }\bar{n}+\left( \frac{\delta }{%
\kappa }\right)^{\prime }\bar{e}
$$%
is a parametrization of the striction curve of $\Phi^{\ast}$.

\section{The relative image of an asymptotically normalized ruled surface}

In this paragraph we consider a non-conoidal ruled surface $\Phi $,
which is asymptotically normalized by $\bar{y} $ via the support function $%
q=fw^{-1} $. The parametrization (\ref{20}) of $\bar{y} $ shows, that the
asymptotic image $\Psi _{1} $ of $\Phi $ is also a ruled surface, whose
generators are parallel to those of $\Phi $. Then, by a straightforward
computation we can find the following parametrization of its striction curve%
\begin{equation}
\Gamma _{1}:\bar{s}_{1}=-\left( \frac{f}{\delta }\right)^{\prime }\bar{e}+%
\frac{f}{\delta }\bar{n}.  \label{37}
\end{equation}%
Thus, if we put for convenience $\bar{y}=\bar{y}_{1}, $ we can rewrite the
parametrization (\ref{20}) as%
$$
\Psi _{1}:\bar{y}_{1}=\bar{s}_{1}+v_{1}\,\bar{e},\quad v_{1}:=-H\,v,
$$
where $H $ denotes the relative mean curvature of $\Phi $ (see (\ref{59})).
Obviously $\Psi _{1} $ is parametrized like in (\ref{4}) and (\ref{23}). We
use $\{\bar{e},\bar{n},\bar{z}\}$ as moving frame of $\Psi _{1} $. The
fundamental invariants of $\Psi _{1} $ are given by%
\begin{equation}
\delta _{1}=-\delta H,\quad \kappa _{1}=\kappa ,\quad \lambda _{1}=-\frac{%
\left( \frac{f}{\delta }\right)^{\prime \prime }+\frac{f}{\delta }}{\kappa
\frac{f}{\delta }}.  \label{35}
\end{equation}%
From the above the following results, which can be checked fairly easily are listed:
\vspace{2mm}
\begin{itemize}
            \item \emph {If $\Phi$ and its asymptotic image $\Psi _{1}$ are congruent} ($\delta=\delta _{1},\kappa =\kappa _{1},\lambda =\lambda _{1})$, \emph {then}
                 $$
                 f=\frac{\delta ^{2}}{\kappa }\quad \text{and}\quad \left( \frac{\delta }{\kappa }\right) ^{\prime \prime }+\frac{\delta }{\kappa }\left( 1+\kappa \lambda \right) =0,
                 $$%
                    \emph{and thus} $\Phi $ \emph {is a proper relative sphere} (see Proposition (\ref{Prop})).          \vspace{1mm}
            \item \emph {$\Psi _{1}$ is orthoid} ($\lambda_{1}=0)$ \emph {iff} $f=\delta \left( c_{1}\cos u+c_{2}\sin u\right),\,c_{1},c_{2}\in \mathbb{R}, \,c_{1}^2+c_{2}^2\neq 0$.\vspace{1mm}
            \item \emph{The striction curve of $\Psi _{1}$ is an asymptotic line of it} ($\kappa _{1}=\lambda _{1}$) \emph {iff}
                $$
                \left( \frac{f}{\delta }\right) ^{\prime \prime }+\frac{f}{\delta }\left(1+\kappa ^{2}\right) =0,
                $$
                and \emph{it is an Euclidean line of curvature of it} ($1+\kappa _{1}\lambda _{1}=0$ ) \emph {iff} $ f=\delta \left( c_{1}u+c_{2}\right)$, \,$c_{1},c_{2}\in \mathbb{R},\,\,c_{1}^2+c_{2}^2\neq0$.\vspace{1mm}
            \item \emph{$\Psi _{1}$ is an Edlinger surface}\footnote{i.e. its osculating quadrics are rotational hyperboloids \cite{Hoschek}}
                ($\delta _{1}^{\prime }=1+\kappa _{1}\lambda _{1}=0$ \cite[p. 36]{Hoschek}) \emph {iff}
                $$
                f=\frac{c\delta }{\kappa }\quad \text{and}\quad \kappa =\frac{1}{c_{1}u+c_{2}
                },\quad c,c_{1},c_{2}\in \mathbb{R},\quad c\neq0,\,\,\,c_{1}^2+c_{2}^2\neq 0.
                $$
\end{itemize}
For $f=|\delta |^{1/2}$, i.e. for the equiaffine normalization, some of the
above results were obtained in \cite[\S \ 4]{Stamou2012}.

\vspace{2mm}
We now assume that $\Phi$ has a ``precedent'' ruled surface, i.e. that there exists another skew ruled surface, say $\Psi ^{\ast }$, with parallel generators, an asymptotic image of which is $\Phi $. We consider a parametrization of $\Psi ^{\ast } $ like in (\ref{4})--(\ref{23}) and let $\delta ^{\ast },\kappa ^{\ast },\lambda ^{\ast} $ be its fundamental invariants. We denote likewise all magnitudes of $\Psi ^{\ast } $ by the usual symbols supplied with a star ($^{\ast }$). We normalize $\Psi ^{\ast } $ asymptotically via the support function $q^{\ast}=f^{\ast }w^{\ast ^{-1}} $, and suppose that the resulting normalization of it, say $\Psi ^{\ast \ast }, $ is the given ruled surface $\Phi $. Then, on account of (\ref{35}), clearly $\kappa ^{\ast }=\kappa $ and
\begin{equation}
\delta =-\delta ^{\ast }H^{\ast },\quad \lambda =-\frac{\left( \frac{f^{\ast
}}{\delta ^{\ast }}\right)^{\prime \prime }+\frac{f^{\ast }}{\delta ^{\ast }}%
}{\kappa \frac{f^{\ast }}{\delta ^{\ast }}},  \label{29}
\end{equation}
where, in view of (\ref{59}), $H^{\ast }=-\delta ^{\ast ^{-2}}\,\kappa \,f^{\ast}$ is the relative mean curvature of $\Phi^{\ast}$. Thus the system (\ref{29}) becomes
\begin{equation}
\frac{f^{\ast }}{\delta ^{\ast }}=\frac{\delta }{\kappa },\quad \left( \frac{%
\delta }{\kappa }\right)^{\prime \prime }+\frac{\delta }{\kappa }\left(
1+\kappa \lambda \right) =0.  \label{30}
\end{equation}%
Let, conversely, the relations (\ref{30}) be valid. We consider an arbitrary skew ruled surface $\Psi^{\ast}$, whose generators are parallel to those of $\Phi$, and let $\delta^{\ast}$ be its distribution parameter. The conical curvature of $\Psi^{\ast}$ equals $\kappa$. We normalize asymptotically $\Psi^{\ast}$ via the support function $q^{\ast}=f^{\ast}\,w^{\ast^{-1}}$, where $f^{\ast}=\delta\,\delta^{\ast}\,\kappa^{-1}$. We can easily verify, by using (\ref{35}) and (\ref{30}), that the fundamental invariants of the asymptotic image $\Psi^{\ast\ast}$ of $\Psi^{\ast}$ coincide with the corresponding fundamental invariants of $\Phi$. Hence $\Psi^{\ast\ast}$ and $\Phi$ are congruent. So we arrive at

\begin{proposition}
The ruled surface $\Phi$ is the asymptotic image of a ruled surface $\Psi ^{\ast }$ iff the second of the conditions (\ref{30}) is
valid.
\end{proposition}

We suppose now that $\Phi $ is not a proper relative sphere ($\Phi \neq
\Psi _{1}) $ and we normalize asymptotically its asymptotic image $\Psi _{1}
$. Let $q_{1}=f_{1}w_{1}^{-1} $ be the support function of $\bar{y}_{1} $.
Analogously to the computations above we get the following parametrization
of the asymptotic image $\Psi _{2} $ of $\Psi _{1}: $%
$$
\Psi _{2}:\bar{y}_{2}=\bar{s}_{2}+v_{2}\,\bar{e},\quad
v_{2}:=-H_{1}\,v_{1},\quad H_{1}=\frac{f_{1}}{fH},
$$%
where%
$$
\Gamma _{2}:\bar{s}_{2}=-\left( \frac{f_{1}}{\delta _{1}}\right)^{\prime }%
\bar{e}+\frac{f_{1}}{\delta _{1}}\bar{n}
$$%
is its striction curve and $H_{1} $ is the relative mean curvature of $\Psi
_{1} $. Thus $\Psi _{2} $ is parametrized like in (\ref{4}) and (\ref{23}).
Obviously, the Tchebychev vector $\bar{T}_{1} $ of $\bar{y}_{1} $ is parallel to $%
\bar{e} $. The fundamental invariants of $\Psi _{2} $ are computed by (see (\ref{35}))
$$
\delta _{2}=-\delta _{1}H_{1},\quad \kappa _{2}=\kappa ,\quad \lambda _{2}=-%
\frac{\left( \frac{f_{1}}{\delta _{1}}\right)^{\prime \prime }+\frac{f_{1}}{%
\delta _{1}}}{\kappa \frac{f_{1}}{\delta _{1}}}.
$$%
According to Proposition (\ref{Prop}) we have: \emph{The asymptotic image $\Psi_{1}$
of $\Phi$ is a proper relative sphere iff there exists a constant $c\neq 0$,
such that $cf_{1}=fH$} (the condition (\ref{32}) is identically fulfilled). Thus, we obtain the following results:\vspace{1mm}
\begin{itemize}
        \item \emph{$\Phi $ and $\Psi _{2}$ are congruent iff}
        $$
        f_{1}=f\quad \text{and}\quad \left( \frac{\delta }{\kappa }\right) ^{\prime
        \prime }+\frac{\delta }{\kappa }\left( 1+\kappa \lambda \right) =0.
        $$
        \item \emph{$\Psi _{1}$ and $\Psi _{2}$ are congruent iff} $\delta ^{2}f_{1}=\kappa f^{2}$.\vspace{2mm}
        \item \emph{$\Psi _{2}$ is orthoid iff} $f_{1}=\frac{\kappa f}{\delta }\left( c_{1}\cos u+c_{2}\sin u\right)
        ,\,c_{1},c_{2}\in \mathbb{R},\,\,c_{1}^2+c_{2}^2\neq 0$.\vspace{2mm}
        \item \emph{The stiction curve of $\Psi _{2}$ is an asymptotic line of it iff}
        $$
        \left( \frac{\delta f_{1}}{\kappa f}\right) ^{\prime \prime }+\frac{\delta
        f_{1}}{\kappa f}\left( \kappa ^{2}+1\right) =0,
        $$
        and \emph {it is an Euclidean line of curvature of it iff}
        $$
        f_{1}=\frac{\kappa f}{\delta }\left( c_{1}u+c_{2}\right) ,\quad c_{1},c_{2}\in \mathbb{R},\quad c_{1}^2+c_{2}^2\neq 0.
        $$
        \item \emph{$\Psi _{2}$ is an Edlinger surface iff}
        $$
        f_{1}=\frac{cf}{\delta }\quad \text{and}\quad \kappa =\frac{1}{c_{1}u+c_{2}} ,\quad c,c_{1},c_{2}\in \mathbb{R}, \quad c\neq 0, \quad c_{1}^2+c_{2}^2\neq 0.
        $$
\end{itemize}

Continuing in the same way we obtain a sequence $\left \{ \Psi _{i}\right \}
_{i\in \mathbb{N}}$ of ruled surfaces, such that $\Psi _{i}$ is the
asymptotic image of $\Psi _{i-1}$. Moreover, if $q_{i-1}=f_{i-1}\,w_{i-1}^{-1}$
is the asymptotic support function of $\Psi _{i-1}$, we can easily check
that the parametrization of $\Psi _{i}$ reads%
$$ \Psi _{i}:\bar{y}_{i}=\bar{s}_{i}+v_{i}\ \bar{e},\quad v_{i}:=-H_{i-1}\,v_{i-1},$$ where
$$
\Gamma _{i}:\bar{s}_{i}=-\left( \frac{f_{i-1}}{\delta _{i-1}}\right)
^{\prime }\bar{e}+\frac{f_{i-1}}{\delta _{i-1}}\bar{n}
$$%
is its striction curve and $H_{i-1}$ is the relative mean curvature of $\Psi
_{i-1}$. $\Psi _{i}$ is parametrized like in (\ref{4}) and (\ref{23}) and
its fundamental invariants are computed by
$$
\delta _{i}=-\delta _{i-1}H_{i-1},\quad \kappa _{i}=\kappa ,\quad \lambda
_{i}=-\frac{\left( \frac{f_{i-1}}{\delta _{i-1}}\right) ^{\prime \prime }+%
\frac{f_{i-1}}{\delta _{i-1}}}{\kappa \frac{f_{i-1}}{\delta _{i-1}}}.
$$%
The relative magnitudes of $\Psi _{i-1}$ are recursively computed by%
$$
J_{i-1}=0,\quad H_{i-1}=S_{i-1}=\frac{f_{i-1}}{f_{i-2}H_{i-2}},\quad K_{i-1}=H_{i-1}^2.
$$%
Finally, we notice that the Tchebychev vectors of all asymptotic normalizations of the sequence $\left \{ \Psi _{i}\right \}
_{i\in \mathbb{N}}$ are parallel to $\bar{e}$ and that their asymptotic developables
 coincide with the director cone of $\Phi$ \cite[p. 263]{Pottmann}.

\section{Some results on the Tchebychev and the support vector fields}

We consider a ruled surface $\Phi$, which is asymptotically normalized by $\bar{y}$
via the support function $q=fw^{-1}$. The Tchebychev vector of $\bar{y}$
can be computed by using (\ref{13}) and (\ref{19}). We find%
$$
\bar{T}=\frac{2\delta f^{\prime }-\delta ^{\prime }f}{2\delta ^{2}}\bar{e}.
$$%
The divergence $\divz ^{I}\bar{T} $ and the rotation $\curl^{I}%
\bar{T} $ of $\bar{T} $ with respect to the first fundamental form $I $ of $%
\Phi$, which initially read \cite[p. 304, 305]{Stamou2012}
$$
\divz ^{I}\bar{T}=\frac{\left( wT^{i}\right) _{/i}}{w},\quad \curl%
^{I}\bar{T}=\frac{\left( g_{12}T^{1}+g_{22}T^{2}\right) _{/1}-\left(
g_{11}T^{1}+g_{12}T^{2}\right) _{/2}}{w},
$$%
become (see (\ref{47}) and (\ref{43}))%
$$
\divz ^{I}\bar{T}=\frac{v\left( 2\delta f^{\prime }-\delta ^{\prime
}f\right) }{2\delta ^{2}w^{2}},\quad \curl^{I}\bar{T}=\frac{\delta
\left( 2\delta f^{\prime \prime }-3\delta^{\prime }f^{\prime }\right)
+f\left( 2\delta ^{\prime 2}-\delta \delta ^{\prime \prime }\right) }{%
2\delta ^{3}w},
$$%
from which we obtain:
\begin{itemize}
        \item \emph{It is } $\divz ^{I}\bar{T}\equiv0$ \emph{iff} $f=c|\delta |^{1/2},\, c\in \mathbb{R}^{\ast}$,
        or equivalently \emph{iff $\bar{T}=\bar{0}$}.\vspace{1mm}
        \item \emph{It is $\curl^{I}\bar{T}\equiv 0$ iff $\delta \left( 2\delta
        f^{\prime \prime }-3\delta ^{\prime }f^{\prime }\right) +f\left( 2\delta
        ^{\prime 2}-\delta \delta ^{\prime \prime }\right) =0$}, or, after standard
        calculation, \emph{iff $f=|\delta |^{1/2}\left( c_{1}\int |\delta
        |^{1/2}du+c_{2}\right)$, $c_{1},c_{2}\in \mathbb{R},\,c_{1}^2+c_{2}^2\neq 0$}.
\end{itemize}

\vspace{2mm}
Let $\divz ^{G}\bar{T} $ and $\curl^{G}\bar{T} $ be the divergence
and the rotation of $\bar{T} $ with respect to the relative metric. In
analogy to the computation above we get
$$
\divz ^{G}\bar{T}\equiv 0,\quad \curl^{G}\bar{T}\equiv 0.
$$%
The relation $\curl^{G}\bar{T}\equiv 0 $ agrees with $\bar{T}=\nabla
^{G}\left( f|\delta |^{-1/2},\bar{x}\right) $ (see (\ref{61})).

\vspace{2mm}
The support vector $\bar{Q} $ of an asymptotic normalization becomes (see (\ref{15}))
\begin{equation}
\bar{Q}=w\frac{\kappa fv+\delta^{\prime }f-\delta f^{\prime }}{4\delta ^{2}f}%
\bar{e}+\frac{v}{4\delta w}\left( v\bar{n}+\delta \bar{z}\right) .
\label{41}
\end{equation}%
We observe, that $\langle \bar{e},\bar{Q}\rangle =0 $ iff
$$
\kappa fv+\delta ^{\prime }f-\delta f^{\prime }=0.
$$
On differentiating twice relative to $v $ we obtain the system
$$
\kappa f=\delta^{\prime }f-\delta f^{\prime }=0,
$$
which implies $\kappa =0$ and $f=c|\delta|,\ c\in \mathbb{R}^{\ast}$. The
inverse also holds. So we have: \emph{The support vectors $\bar{Q}$ are orthogonal to the generators iff $\Phi$ is
conoidal and $f=c\,|\delta|$, $c\in \mathbb{R}^{\ast}$.}
On account of (\ref{45}) a direct computation yields%
\begin{equation}
\divz ^{I}\bar{Q}=\frac{3\kappa fv^{2}+\left( \delta ^{\prime }f-2\delta
f^{\prime }\right) v+\delta ^{2}f\left( \kappa -\lambda \right) }{4\delta
^{2}fw},  \label{53}
\end{equation}%
\begin{equation}
\curl^{I}\bar{Q}=\frac{A_{3}v^{3}+A_{2}v^{2}+A_{1}v+A_{0}}{4\delta
^{3}f^{2}w^{2}},  \label{58}
\end{equation}%
where
\begin{subequations}
\begin{eqnarray}
A_{3}&=&f^{2}\left(\delta \kappa ^{\prime }-2\delta ^{\prime }\kappa
\right) , \\
A_{2} &=&-2\delta^{\prime 2}f^{2}+\delta f \left(\delta^{\prime} f^{\prime}+\delta ^{\prime \prime }f\right) +\delta^{2}\left[f^{\prime 2}-2f^{2}\left(1+\kappa \lambda \right)  -ff^{\prime \prime }\right], \\
A_{1} &=&\delta ^{2}f \big[ \delta \lambda f^{\prime }+f\left[ \delta
\kappa^{\prime }-\delta ^{\prime }\left( \kappa +\lambda \right) \right] \big] , \\
A_{0} &=&-\delta^{2}\big[ f^2(\delta^{\prime 2}-\delta\delta^{\prime \prime})+\delta^{2} [ff^{\prime \prime}+f^2\left( 1+\kappa \lambda \right)-f^{\prime 2} ] \big] .
\end{eqnarray}%
\end{subequations}
Also we have
\begin{equation}
\divz ^{G}\bar{Q}=\frac{2\kappa fv^{4}+\left( \delta ^{\prime }f-2\delta
f^{\prime }\right) v^{3}+3\delta ^{2}\kappa fv^{2}-2\delta^{3}f^{\prime}v+\delta^{4}f \left(\kappa -\lambda \right) }{4\delta ^{2}fw^{3}},  \label{54}
\end{equation}%
and%
\begin{equation}
\curl^{G}\bar{Q}\equiv 0.  \label{55}
\end{equation}%
Let $\divz ^{I}\bar{Q}=0 $. Then by (\ref{53}) we have
$$3\kappa
fv^{2}+\left( \delta ^{\prime }f-2\delta f ^{\prime }\right) v+\delta
^{2}f\left( \kappa -\lambda \right) =0, $$
from which, by successive differentiations relative to $v, $ we infer the system
 $$
\kappa f =\delta ^{\prime }f - 2\delta f^{\prime }= \delta^{2} f \left(\kappa -\lambda \right) =0,
$$
i.e. $\kappa =\lambda =0$ and $f=c\,|\delta|^{1/2},\,c\in\mathbb{R}^{\ast}$.
The inverse also holds. So we have: \emph{It is} $\divz^ {I}\bar{Q}\equiv0$ \emph{iff $\Phi$ is a right conoid and $f=c\,|\delta|^{1/2},c \in\mathbb{R}^{\ast} $.}
Treating the relations (\ref{58})--(\ref{55}) similarly we obtain the following results:\vspace{1mm}
\begin{itemize}
                \item \emph {It is }$\curl^{I}\bar{Q}\equiv 0$ \emph{iff}
                    \begin{itemize}
                    \item $\Phi$ \emph{is an Edlinger surface with constant invariants and $f=c\in \mathbb{R}^{\ast }$}, or
                    \item $\Phi$ \emph{is a right conoid}, $\delta=\frac{c_{1}}{u+c_{2}}$ \emph{and} $f=\frac{c_{1}c_{3}}{(u+c_{2}) \sqrt{e^{u(u+2c_{2})}}}$, $c_{1},c_{3}\in \mathbb{R}^{\ast},\,c_{2}\in \mathbb{R}$,  or
                    \item \emph{the fundamental invariants of $\Phi$ fulfil the relations}\\
                    $$
                     c_{1}^{2}\delta ^{6}-5c_{3}\left [ \delta \left( u + c_{1}\right) + c_{3}\right]=0,\quad
                    \kappa =c_{1}\delta ^{2},\quad \lambda =\frac{-c_{1}\delta ^{4}}{c_{3}^{2}+c_{1}^{2}\delta ^{6}},\quad c_{1},c_{2},c_{3}\in \mathbb{R}^{\ast},
                    $$%
                    \emph{and} $f=c_{2}\,|\delta |\,e^{c_{3}\int \frac{du}{\delta }}$.
                    \end{itemize}
                \item \emph{It is }$\divz ^{G}\bar{Q}\equiv 0$ \emph{iff $\Phi$ is a right helicoid and $f=c\in \mathbb{R}^{\ast}$.}
\end{itemize}

\vspace{2mm}

We consider now the following families of curves on $\Phi $:
a) the curved asymptotic lines,
b) the curves of constant striction distance $(u $-curves) and
c) the $\widetilde{K} $-curves, i.e. the curves along which the Gaussian
curvature is constant \cite{Sachs}.
The corresponding differential equations of these families of curves are%
\begin{equation}
\kappa v^{2}+\delta'v+\delta^{2}\left( \kappa -\lambda \right) -2\delta
v^{\prime }=0,  \label{38}
\end{equation}%
\begin{equation}
v^{\prime }=0,  \label{39}
\end{equation}%
\begin{equation}
2\delta vv^{\prime }+\delta ^{\prime }\left( \delta ^{2}-v^{2}\right) =0.
\label{40}
\end{equation}%
It will be our task to investigate necessary and sufficient conditions for
the support vector field $\bar{Q}$ to be tangential or orthogonal to one of
these families of curves. To this end we consider a directrix $\Lambda
:v=v\left( u\right) $ of $\Phi $. Then we have
\begin{equation}
\bar{x}^{\prime }=\left( \delta \lambda +v^{\prime }\right) \bar{e}+v\bar{n}%
+\delta \bar{z}.  \label{36}
\end{equation}%
From (\ref{41}) and (\ref{36}) it follows: $\bar{x}^{\prime }$ and $\bar{Q}$
are parallel or orthogonal iff%
\begin{equation}
\kappa fv^{3}+\left( \delta ^{\prime }f-\delta f^{\prime }\right)
v^{2}+\delta f\left[ \delta \left( \kappa -\lambda \right) -v^{\prime }%
\right] v+\delta ^{2}\left( \delta ^{\prime }f-\delta f^{\prime }\right) =0
\label{42}
\end{equation}%
or%
\begin{equation}
\left( \kappa fv+\delta ^{\prime }f-\delta f^{\prime }\right) \left( \delta
\lambda +v^{\prime }\right) +\delta fv=0,  \label{56}
\end{equation}%
respectively. Then, from (\ref{38}) and (\ref{42}) (resp. (\ref{56})), we
infer, that $\bar{Q}$ is tangential or orthogonal to the curved asymptotic
lines iff%
\begin{equation}
\kappa fv^{3}+\left( \delta ^{\prime }f-2\delta f^{\prime }\right)
v^{2}+\delta ^{2}f\left( \kappa -\lambda \right) v+2\delta ^{2}\left( \delta
^{\prime }f-\delta f^{\prime }\right) =0  \label{57}
\end{equation}%
or%
\begin{equation}
\kappa ^{2}fv^{3} {\small +} \kappa \left( 2\delta ^{\prime }f {\small -}\delta f^{\prime
}\right) v^{2} {\small +}[\delta ^{2}\kappa f\left( \kappa {\small +}\lambda \right) {\small +}\delta
^{\prime }\left( \delta ^{\prime }f{\small -}\delta f^{\prime }\right) {\small +}2\delta
^{2}f]v {\small +}\delta ^{2}\left( \delta ^{\prime }f{\small -}\delta f^{\prime }\right)
\left( \kappa {\small +}\lambda \right) =0,  \label{62}
\end{equation}%
respectively. From (\ref{57}) and (\ref{62}), after successive differentiations relative to $v,$ we obtain%
$$
\kappa f =\delta ^{\prime }f-2\delta f^{\prime }=\delta^{2} f \left(\kappa -\lambda\right) = 2\delta^{2} \left(\delta
^{\prime }f-\delta f^{\prime }\right)=0
$$%
and%
$$
\kappa^{2} f =\kappa \left( 2\delta ^{\prime }f -\delta f^{\prime }\right) =\delta
^{2}\kappa f\left( \kappa + \lambda \right) +\delta ^{\prime }\left( \delta
^{\prime }f -\delta f^{\prime }\right) +2\delta ^{2}f = \delta^{2} \left( \delta ^{\prime
}f -\delta f^{\prime }\right) \left( \kappa +\lambda \right) =0,
$$%
respectively. Standard treatment of these systems leads to the following results:
\vspace{2mm}
\begin{itemize}
  \item \emph {$\bar{Q}$ is tangential to the curved asymptotic lines of $\Phi$ iff $\Phi$
is a right helicoid and $f {\small{=}}c\in \mathbb{R}^{\ast }$}.\vspace{1mm}
  \item \emph{$\bar{Q}$ is orthogonal to the curved asymptotic lines of $\Phi $ iff $%
\Phi $ is a right conoid and the function $f$ is given by $f=c\,|\delta|\,e^{2\int \frac{\delta }{\delta ^{\prime }}du},c\in \mathbb{R}^{\ast }$}.
\end{itemize}

\vspace{2mm}
From (\ref{39}) and (\ref{42}), resp. (\ref{56}), we obtain: $\bar{Q} $ is
tangential or orthogonal to the $u $-curves iff%
$$
\kappa fv^{3}+\left( \delta ^{\prime }f-\delta f^{\prime }\right)
v^{2}+\delta ^{2}f\left( \kappa -\lambda \right) v+\delta ^{2}\left( \delta
^{\prime }f-\delta f^{\prime }\right) =0
$$%
or%
$$
f\left( 1+\kappa \lambda \right) v+\lambda \left( \delta ^{\prime }f-\delta
f^{\prime}\right)=0,
$$%
respectively. Treating these polynomials in the same way we result:
\vspace{2mm}
\begin{itemize}
  \item \emph{$\bar{Q}$ is tangential to the $u$-curves of $\Phi$ iff $\Phi$ is a
right conoid and $f=c\,|\delta |,c\in \mathbb{R}^{\ast }$}.\vspace{1mm}
  \item \emph{$\bar{Q}$ is orthogonal to the $u$-curves of $\Phi $ iff the striction
curve of $\Phi $ is an Euclidean line of curvature and $f=c\,|\delta |,c\in
\mathbb{R}^{\ast }$}.
\end{itemize}

\vspace{1mm}
From (\ref{40}) and (\ref{42}), resp. (\ref{56}), we obtain: $\bar{Q}$ is
tangential or orthogonal to the $\widetilde{K}$-curves iff%
$$
2\kappa fv^{3} +\left(\delta ^{\prime }f-2\delta f^{\prime }\right)
v^{2}+2\delta ^{2}f\left( \kappa -\lambda \right) v +\delta ^{2}\left(
3\delta ^{\prime }f-2\delta f^{\prime }\right)=0
$$
or
$$
{\delta ^{\prime }\kappa fv^{3}{\small +}\left[2\delta^{2}f\left(1{\small +}\kappa \lambda
\right) {\small +}\delta ^{\prime }\left( \delta ^{\prime }f {\small -}\delta
f^{\prime}\right) \right] v^{2}}
{\small +}\delta ^{2}\left[ \delta ^{\prime }f\left( 2\lambda {\small -}\kappa \right)
{\small -}2\delta \lambda f^{\prime }\right] v{\small -}\delta ^{2}\delta ^{\prime }\left(
\delta ^{\prime }f{\small -}\delta f^{\prime }\right)=0,
$$
respectively. Treating analogously these polynomials we easily obtain:
\vspace{2mm}
\begin{itemize}
  \item  \emph{$\bar{Q}$ is tangential to the $\widetilde{K}$-curves of $\Phi $ iff $\Phi$ is a right helicoid and $f=c\in \mathbb{R}^{\ast }$}.\vspace{1mm}
  \item  \emph{$\bar{Q}$ is orthogonal to the $\widetilde{K}$-curves of $\Phi $ iff $\Phi$ is an Edlinger surface and $f=c\in \mathbb{R}^{\ast}$}.
\end{itemize}

\vspace{2mm}
To complete this work we consider the Euclidean lines of curvature of $\Phi$. Their differential equation, initially being%
$$
g_{12}h_{11}-g_{11}h_{12}+\left( g_{22}h_{11}-g_{11}h_{22}\right) v^{\prime
}+\left( g_{22}h_{12}-g_{12}h_{22}\right) v^{\prime 2}=0,
$$%
becomes, on account of (\ref{47}) and (\ref{7}),%
$$
\delta \left[ w^{2}\left( 1+\kappa \lambda \right) +\delta ^{\prime }\lambda
v\right] +\left[ \kappa w^{2}+\delta ^{\prime }v-\delta ^{2}\lambda \right]
v^{\prime }-\delta v^{\prime 2}=0,
$$%
from which, by virtue of (\ref{42}), we infer, that $\bar{Q}$ is tangent
to the one family of the lines of curvature of $\Phi $ iff
$$
{{\small -}\kappa ff^{\prime }v^{3}{\small +}\left[ \delta f^{\prime 2}{\small -}\delta f^{2}\left(
1{\small +}\kappa \lambda \right) {\small -}\delta ^{\prime }ff^{\prime }\right] v^{2}} {\small +}\delta f\left( \kappa {\small -}\lambda \right) \left( \delta ^{\prime }f{\small -}\delta
f^{\prime }\right) v{\small +}\delta \left( \delta f^{\prime }{\small -}\delta ^{\prime
}f\right) ^{2}=0.
$$
It results the system%
$$
\kappa f f^{\prime }=
\left[ \delta f^{\prime 2}{\small -}\delta f^{2}\left(1{\small +}\kappa \lambda \right) {\small -}\delta ^{\prime }ff^{\prime }\right]=
\delta f\left( \kappa {\small -}\lambda \right) \left( \delta ^{\prime }f{\small -}\delta f^{\prime }\right) =
 \delta \left( \delta f^{\prime }{\small -}\delta ^{\prime}f\right) ^{2}=0,
$$
from which we get
$$
\delta ^{\prime }=1+\kappa \lambda =f^{\prime }=0.
$$%
Hence $\Phi $ \emph{is an Edlinger surface and the function $f$ is constant}. Moreover, we can easily confirm, that the Euclidean principal directions at a point $P$ of an Edlinger surface read%
$$
v^{\prime }=0\quad \text{and}\quad v^{\prime }=\frac{\delta ^{2}+\kappa^{2}w^{2}}{\delta \kappa }.
$$%
Since the second of these relations verifies (\ref{42}), we have: \emph{When} $\Phi$ \emph{is an Edlinger surface and the function} $f$ \emph{is constant, then the support vector field} $\bar{Q}$ \emph{is tangent to those Euclidean lines of
curvature of }$\Phi $, \emph{which are orthogonal to the striction curve of} $\Phi$.

\end{document}